\documentclass[12pt,a4paper]{amsart}
\usepackage{amssymb}
\makeatletter
\renewcommand{\@begintheorem}[2]{
\rm \trivlist \item [\hskip \labelsep {\bf #2\ \ #1.}]
                                }
\makeatother

\makeatletter

\DeclareFontFamily{U}{cyr}{}
\DeclareFontShape{U}{cyr}{m}{n}{
  <5> wncyr5 <6> wncyr6 <7> wncyr7 <8> wncyr8 <9> wncyr9 <10->
wncyr10}{}
\DeclareMathAlphabet{\mathcyr}{U}{cyr}{m}{n}

\input cyracc.def
\input epsf

\newcommand{\ts}{\vspace{\baselineskip}\noindent{\bf Proof.}$\;\;$}
\newcommand{\ZZ}{{\bf Z}}
\newcommand{\QQ}{{\bf Q}}

\newcommand{\CC}{{\bf C}}

\newcommand{\PP}{{\bf P}}

\newcommand{\bes}{\begin{equation*}}
\newcommand{\ees}{\end{equation*}}

\title{Mirror Quintics, discrete symmetries and Shioda Maps}

\author{Gilberto Bini}
\address{Dipartimento di Matematica, Universit\`a di Milano,
Via Saldini 50, I-20133 Milano, Italia}
\email{gilberto.bini@mat.unimi.it}
\author{Bert van Geemen}
\address{Dipartimento di Matematica, Universit\`a di Milano,
Via Saldini 50, I-20133 Milano, Italia}
\email{geemen@mat.unimi.it}
\author{Tyler L.\  Kelly}
\address{Department of Mathematics,
University of Georgia, 325B Boyd, Athens, GA 30605, USA}
\email{tlkelly@uga.edu}

\begin{document}

\begin{abstract}
In a recent paper, Doran, Greene and Judes considered one parameter families of quintic threefolds with finite symmetry groups.
A surprising result was that each of these six families has the same Picard Fuchs equation associated to the holomorphic $3$-form.
In this paper we give an easy argument, involving the family of Mirror
Quintics, which implies this result. Using a construction due to Shioda, we also relate certain quotients of these one parameter families to the family of Mirror Quintics. Our constructions generalize to degree $n$ Calabi Yau varieties in $(n-1)$-dimensional projective space.

\end{abstract}

\maketitle

\section*{Introduction}
Quintic threefolds in projective $4$-space with a finite automorphism
group have been studied for applications to string theory. In particular,
any smooth quintic $X_t$ in the Dwork pencil
(see section \ref{mirrors}) has a
group $H\cong (\ZZ/5\ZZ)^3$ of automorphisms  which act trivially on the
holomorphic three form.
The quotient variety $X_t/H$ has a resolution of singularities $M_t$ which is again a Calabi Yau (CY) threefold and its Hodge numbers are the `mirrors' of those of $X_t$: $h^{p,q}(M_t)=h^{3-p,q}(X_t)$. This was discovered by
Greene and Plesser \cite{GP} and started an ongoing,
extensive and profound study of CY threefolds.

An important ingredient in String theory is the Picard Fuchs equation
satisfied by the periods of the holomorphic $3$-form of the Mirrors. As
these Mirrors are quotients of the quintics in the Dwork pencil, the Picard Fuchs equation for that pencil is the same. In \cite{DGJ}, the Picard Fuchs equations of five other pencils, $X_{A,t}$, of quintic threefolds were determined, where $A$ are certain $5\times 5$ matrices. Somewhat surprisingly, these turned out to be the same as the one for the Mirror family. Here we show that there are maps from $X_{A,t}$ to the Mirror quintic $M_t$; moreover, the Mirror quintics are quotients
of the $X_{A,t}$ by finite groups. Thus the equality of the Picard Fuchs equations is obvious. To show that for each pencil the Mirror quintics are quotients, we follow a construction due to Shioda which gives a pencil
$X_{dI,t}$ of degree $d$ (where $d$ depends on $A$ and in general $d>5$) 3-folds in $\PP^4$ which maps to both the pencil under consideration and the Dwork pencil. Using the action of a finite group on this degree $d$ pencil we obtain the desired results. In the diagram below the maps are rational maps.
$$
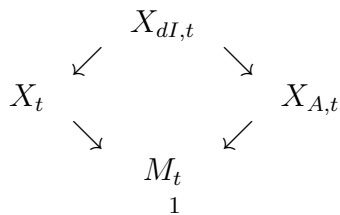
\begin{array}{lcccr}
&&X_{dI,t}&&\\
&\swarrow&&\searrow&\\
X_t&&&&X_{A,t}\\
&\searrow&&\swarrow&\\
&&M_t&&
\end{array}
$$

\section{The Mirror quintics and symmetric quintics}

\subsection{The Mirror quintics}\label{mirrors}
The Mirror quintics are a one parameter family of CY threefolds
$M_t$ which have Hodge numbers $h^{1,1}(M_t)=101$ and $h^{2,1}(M_t)=4$.
Following \cite{GP},
they are defined as (crepant) desingularizations
of the quotients of quintics in the Dwork pencil in $\PP^4$
$$
X_{t}\,:=\,Z(F_{5I,t}),\qquad
F_{t}\,=\,x_1^5+x_2^5+\ldots+x_5^5\,-\,5tx_1x_2\cdots x_5
$$
by the finite group
$$
H\,=\,
\langle\,
h_1:=h_{(-1,1,0,0,0)},h_2=h_{(-1,0,1,0,0)},\,
h_3:=h_{(-1,0,0,1,0)}\,\rangle
\,\cong\,(\ZZ/5\ZZ)^{3}
$$
where, for a primitive fifth root of unity $\zeta=\zeta_5$, we define automorphisms of $\PP^4$ by:
$$
h_{(a_1,\ldots,a_5)}(x_1:\ldots:x_5)\,=\,
(\zeta^{a_1}x_1:\ldots:\zeta^{a_5}x_5)
\qquad (\zeta^5=1,\quad\zeta\neq 1).
$$

The (singular) quotients of the Dwork quintics by $H$ can be found as follows. The generators $h_1,h_2,h_3$ of $H$
act on the affine open subset
$U_{t}:=\{\,(x_1:\ldots:x_4:1)\in X_{t}\,\}$.
The quotient of this affine variety is (by definition)
the spectrum of the subring of $H$-invariants in the affine coordinate ring $\CC[U_t]$:
$$
U_t/H\,:=\,Spec(\CC[U_t]^H),\qquad
\CC[U_t]=\CC[x_1,\ldots,x_4]/(F_{t}(x_1,\ldots,x_4,1)).
$$
As elements of $H$ fix any monomial up to a fifth root of unity, the
ring $\CC[U_t]^H$ is generated by $\CC$ and the invariant monomials.
A monomial $x_1^{k_1}\ldots x_4^{k_4}$ is invariant under $H$ iff
$k_1\equiv k_2\equiv k_3\equiv k_4$ mod $5$, thus an invariant monomial is an element of the subring $\CC[z_0,z_1,\ldots,z_5]$ of $\CC[U_t]$
where $z_0=x_1x_2x_3x_4$ and $z_i:=x_i^5$ for $i=1,\ldots,4$.
These invariants satisfy one linear relation obtained from $F_{t}$: $L_t:=z_1+\ldots+z_4+1-5tz_0=0$, in particular, if $t\neq 0$ the ring of invariants is generated by $z_1,\ldots,z_4$.
For any $t$ we also have $z_1\cdots z_4=(x_1\cdots x_4)^5=z_0^5$
hence
$$
\CC[U_t]^H\,\cong\,\CC[z_0,z_1,\ldots,z_4]/(L_t,G),\qquad
G\,:=\,z_0^5-z_1z_2z_3z_4.
$$
Using the $S_5$-symmetry of the equation defining the Dwork pencil it then follows that the quotient $X_{t}/H$, which we denote by $\overline{M}_t$,
is the following variety:
$$
\overline{M}_t:=X_{t}/H\,\cong\,Z(\,z_1+z_2+\ldots+z_5-5tz_0,\;
z_0^5-z_1z_2\cdots z_5\,)\qquad(\subset\PP^5).
$$
For $t\neq 0$, this can be simplified to $Z((z_1+z_2\ldots+z_5)^5-(5t)^5z_1z_2\cdots z_5)\subset\PP^4$ (with homogeneous coordinates $z_1,\ldots,z_5$). In particular, the Mirrors
have a singular birational model which is again a quintic in $\PP^4$.
The quotient map is given by
$$
q_t:\,X_t\,\longrightarrow\,\overline{M}_t,\qquad
(x_1:\ldots:x_5)\,\longmapsto\,
(z_0:z_1:\ldots:z_5):=(x_1x_2\cdots x_5:x_1^5:\ldots:x_5^5).
$$

\subsection{Symmetric quintics}\label{symmetric}
In the paper \cite{DGJ} the following six one-parameter families of
quintics in $\PP^4$ are considered (the last column will be explained in section \ref{XA}).
{\renewcommand{\arraystretch}{1.3}
$$
\begin{array}{|@{\hspace{6pt}}l||@{\hspace{35pt}}r|@{\hspace{10pt}}r|}
\hline
&
F_{A,t}&
d
\\\hline \hline
1&x_1^5+x_2^5+\ldots+x_5^5-5tx_1x_2\cdots x_5&
5 \\ \hline
2&x_1^4x_2+x_2^4x_3+x_3^4x_4+x_4^4x_5+x_5^4x_1-5tx_1x_2\cdots x_5&
5^2\cdot 41=1025\\ \hline
3&x_1^4x_2+x_2^4x_3+x_3^4x_4+x_4^4x_1+x_5^5-5tx_1x_2\cdots x_5&
5\cdot 3\cdot 17=255 \\ \hline
4&x_1^4x_2+x_2^4x_3+x_3^4x_1+x_4^5+x_5^5-5tx_1x_2\cdots x_5 &
5\cdot 13=65\\ \hline
5&x_1^4x_2+x_2^4x_3+x_3^4x_1+x_4^4x_5+x_5^4x_4-5tx_1x_2\cdots x_5&
5\cdot 3\cdot 13=195\\ \hline
6&x_1^4x_2+x_2^4x_1+x_3^5+x_4^5+x_5^5-5tx_1x_2\cdots x_5&
5\cdot 3=15\\ \hline
\end{array}
$$
}
The equation defining each of these pencils is of the form
$$
F_{A,t}\,:=\,\sum_{i=1}^5 \,\prod_{j=1}^5x_j^{a_{ij}}\,-5tx_1x_2\cdots x_5,
\qquad
a_{i1}+a_{i2}+\ldots+a_{i5}=5
$$
for all $i$ (so the sum of the entries in a row of the matrix $A:=(a_{ij})$ is constant and equal to $5$). Moreover, one easily checks that also the sum of the entries in any column of $A$ is equal to $5$:
$ a_{1j}+a_{2j}+\ldots+a_{5j}=5$ for any $j$.
This implies that if we define, for $i=1,\ldots,5$,
$$
z_i\,:=\,\prod_{j=1}^5x_j^{a_{ij}},\qquad\mbox{then}\quad
z_1z_2\cdots z_5=(x_1x_2\cdots x_5)^5.
$$
Now that suffices to show that there is a non-constant rational map
$$
q_{A,t}:\,X_{A,t}\,\longrightarrow\,\overline{M}_t,\qquad
(x_1:\ldots:x_5)\,\longmapsto\,
(z_0:z_1\ldots:z_5),
$$
where $z_0:=x_1x_2\cdots x_5$. The map $q_t$ from the previous section is just $q_t=q_{5I,t}$ where $I$ is the identity matrix.

In Theorem \ref{quotients} we will show that $q_{A,t}$ is also a (birationally) quotient map
and thus $\overline{M}_t$ is birationally isomorphic to $X_{A,t}/H_A$ for a finite group $H_A$ of automorphisms
of $X_{A,t}$.
This result was already suggested in \cite{GPR} for the second family and in \cite{DGJ} for all six families.

\subsection{Picard Fuchs equations}\label{PF}
Of particular interest for applications to string theory are the
Picard Fuchs equations for the periods of the holomorphic 3-forms
on the quintics in these families.
In \cite{DGJ}, an elegant method to find this differential equation
is given. The result is that for all six families one finds the same Picard Fuchs equation, which is also the Picard Fuchs equation of the Mirror Quintics.

The map $q_{A,t}$ is dominant, so its image contains an open subset of $\overline{M}_t$, because it is a``quotient (by a finite group) map"
(see Theorem \ref{quotients}). As $M_t$ is a desingularization of $\overline{M}_t$,
there is then also a dominant rational map $X_{A,t}\rightarrow M_t$.
Thus the holomorphic 3-form on $M_t$ pulls back to a holomorphic 3-form on $X_{A,t}$ and this implies that the Picard Fuchs equations are indeed the same. Actually the map $q_{A,t}$ defines a correspondence between $X_{A,t}$ and the Mirror Quintic $M_t$ which induces an isomorphism of rational Hodge structures between $H^3(M_t,\QQ)$
and a Hodge substructure of $H^3(X_{A,t},\QQ)$. In particular, the variations of these Hodge structures in both one parameter families are isomorphic and hence the Picard Fuchs equations are the same.

\section{The Shioda map}\label{shioda map}
\subsection{The variety $X_A$} \label{XA}
In order to investigate the maps $q_{A,t}$ we introduce a slight generalization of a construction of T.\ Shioda.

Let $A$ be an $n\times n$ matrix with non-negative integer entries
such that the sum of the entries in each row is fixed integer $e$,
which is independent of the row:
$$
A=(a_{ij})\quad(\in M_n(\ZZ)),\qquad a_{ij}\in\ZZ_{\geq 0},\qquad
\sum_{j=1}^n a_{ij}\,=\,e,\quad i=1,\ldots,n.
$$
For such a matrix we define a homogeneous polynomial $F_A$ of degree $e$ in
$n$ variables, a sum of $n$ monomials, its zero locus is a
(not necessarily smooth or irreducible) projective variety $X_A\subset \PP^{n-1}$:
$$
F_A\,:=\,\sum_{i=1}^n\,\prod_{j=1}^n x_j^{a_{ij}}\;=\,
x_1^{a_{11}}x_2^{a_{12}}\ldots x_n^{a_{1n}}\,+\,
x_1^{a_{21}}x_2^{a_{22}}\ldots x_n^{a_{2n}}\,+\,\ldots,\qquad
X_A\,:=\,Z(F_A).
$$
In the case that $A=eI$, we obtain the Fermat hypersurface $X_e$ of degree $e$: $X_{eI}=X_e$.

Now we assume that the matrix $A$ is invertible.
Let $d\in\ZZ_{>0}$ be the smallest positive integer such that $dA^{-1}$ has integer coefficients.
As the cofactor matrix $A^*$ of $A$ has integral entries and
$A^*=\det(A)A^{-1}$, the integer $d$ divides $\det(A)$. We define an $n\times n$ matrix with integral entries $B$ by:
$$
B\,:=\,dA^{-1}\qquad(\in M_n(\ZZ)),\qquad AB=dI=BA.
$$

Consider the column vector ${\bf e}=(1,\ldots,1)\in\ZZ^n$.
The condition $\sum_{j=1}^n a_{ij}\,=\,e$ for all $i$ is then
equivalent to $A{\bf e}=e{\bf e}$.
As $B=dA^{-1}$, this implies that $B{\bf e}=(d/e){\bf e}$,
so the sum of the entries in any row of $B$ is $\sum_kb_{jk}=d/e=:m$.
So any row sum of $B$ is $m$ and $d=me$ is an integer multiple of $e$.

In \ref{symmetric} we already considered the condition that the column sums of $A$ are equal to $n$. As the sum of all entries of $A$ is $n^2$,
the row sums will then be $e=n$.
Having column sums equal to $n$ is equivalent to
${}^t{\bf e}A=n{\bf e}$.
As $A=dB^{-1}$ this implies that ${}^t{\bf e}B=m{\bf e}$, so the column sums of $B$ are constant as well and are equal to $m=d/n$.

\subsection{The Shioda map}\label{shiodamap}
Shioda (\cite{shioda}, p.421) found a rational map $\phi_A$ from the
Fermat variety $X_d\subset \PP^{n-1}$ of degree $d$ to $X_A$ defined by:
$$
\phi_A\,:\,X_d\,\longrightarrow\,X_A,\qquad
(y_1:\ldots:y_n)\,\longmapsto (x_1:\ldots:x_n),\quad
x_j=\prod_{k=1}^ny_k^{b_{jk}}.
$$
To see that this works, substitute the expressions for the $x_j$ in $F_A$ and notice that you get the equation for the Fermat:
{\renewcommand{\arraystretch}{1.3}
$$
\begin{array}{rcl}
\phi_A^*F_A&=&
\sum_{i=1}^n\,\prod_{j=1}^n \left(\prod_{k=1}^ny_k^{b_{jk}}\right)^{a_{ij}}\\
&=&
\sum_{i=1}^n\,\prod_{k=1}^ny_k^{\sum_ja_{ij}b_{jk}}\\
&=&\sum_{i=1}^n\,y_i^d.
\end{array}
$$
} 

\subsection{The generalized Mirror} \label{gen mirror}
We will be interested in the case that
the degree $e$ of $X_A$ is equal to the number of variables $n$ (Calabi-Yau case).
In particular, we consider the generalized Dwork pencil defined by
$$
X_{A,t}\,:=\,Z(F_{A,t}),\qquad
F_{A,t}\,=\,F_A\,-\,nt\left(\prod_{j=1}^n x_j\right).
$$
In case $A=nI$ this is the Dwork pencil of $(n-2)$-dimensional CY varieties in $\PP^{n-1}$, it was studied extensively in \cite{HSBT}.

As we observed in \ref{symmetric}, in case $n=5$ and the column sums of $A$ are constant, then there are rational maps $X_{A,t}\rightarrow \overline{M}_t=X_{dI,t}/H$.
It is straightforward to generalize the results from
\ref{mirrors}: for any $n\times n$ matrix with non-negative coefficients $A$ whose row and column sums are equal to $n$
there is a rational map
$$
q_{A,t}:X_{A,t}\,\longrightarrow\,\overline{M}_t,\quad
(x_1:\ldots:x_n)\,\longmapsto\,
(z_0:z_1\ldots:z_n):=(\prod_{j=1}^n x_j:
\prod_{j=1}^n x_j^{a_{1j}}:\ldots:\prod_{j=1}^n x_j^{a_{nj}}),
$$
where the $(n-1)$-dimensional variety $\overline{M}_t$ is defined by
$$
\overline{M}_t\,:=\,Z(\,z_1+z_2+\ldots+z_n-ntz_0,\;
z_0^n-z_1z_2\cdots z_n\,)\qquad(\subset\PP^n).
$$

\subsection{The generalized Shioda map} \label{deforming}
In order to generalize the Shioda map, we consider the pencil
of degree $d=nm$ varieties in $\PP^{n-1}$ defined by
$$
X_{dI,t}\,:=\,Z(F_{dI,t}),\qquad
F_{dI,t}\,=\,\sum_{j=1}^n y_j^d\,-\,nt\left(\prod_{j=1}^n y_j\right)^m
\qquad(d=mn).
$$
For matrices $A$ as in \ref{gen mirror} we have $\phi_A^*(F_A)=F_{dI}$,
and, moreover, the column sums $\sum_j b_{jk}$ of $B$ are equal to $m$ so:
$$
\phi_{A,t}^*\left(\prod_{j=1}^n x_j\right)\,=\,
\prod_{j=1}^n\left(\prod_{k=1}^ny_k^{b_{jk}}\right)\,=\,
\prod_{k=1}^ny_k^{\sum_j b_{jk}}\,=\,\prod_{k=1}^ny_k^m.
$$
Thus we have (rational) Shioda maps
$\phi_{A,t}:X_{dI,t}\rightarrow X_{A,t}$
and following sequence of rational maps will help to understand the $q_{A,t}$'s:
$$
X_{dI,t}\;\stackrel{\phi_{A,t}}{\longrightarrow}\;
X_{A,t}\;\stackrel{q_{A,t}}{\longrightarrow}\;
\overline{M}_t.
$$
We will show that if $X_{A,t}$ is irreducible then both maps, and their composition, are quotients by finite abelian groups in Theorem \ref{quotients}.

\subsection{Automorphisms}\label{autos}
Let $\zeta=\zeta_d$ be a generator of the cyclic group of $d$-th roots of unity
$$
\mu_d\,:=\,\{z\in\CC\,:\,z^d=1\,\}\,=\,\{\zeta^a\,:a=0,\ldots,d-1\}.
$$
For $a=(a_1,\ldots,a_n)\in(\ZZ/d\ZZ)^n$ we define an automorphism $g_a$
of $\PP^{n-1}$ by
$$
g_a(y_1:\ldots:y_n)\,:=\,(\zeta^{a_1}y_1:\ldots:\zeta^{a_n}y_n).
$$
Note that $a,b\in (\ZZ/d\ZZ)^n$ define the same automorphism iff
$a-b=(k,\ldots,k)$ for some $k\in\ZZ/d\ZZ$. Let
$$
\Gamma_{d}\,:=\,\{g_a\,:\,a=(a_1,\ldots,a_n),\;
a_1+\ldots+a_n\equiv 0\;\mbox{mod}\;n\,\}/
\langle g_{(1,1,\ldots,1)}\rangle.
$$
Generators of this group are the $g_i$, $i=0,\ldots,n-2$, defined as
$$
\Gamma_{d}\,\cong\,
\langle\, g_0:=g_{(n,0,\ldots,0)},\;
g_1:=g_{(-1,1,0\ldots,0)},\ldots,\,
g_{n-2}:=g_{(-1,\ldots,0,1,0)}\,\rangle
\,\cong \,\mu_{m}\times \mu_d^{n-2},
$$
note that $g_0g_1^n=g_{(0,n,0,\ldots,0)}$, $g_0g_1^ng_2^n=g_{(0,0,n,0,\ldots,0)}$ etc.

For any $t$, the group $\Gamma_{d}$ is a subgroup of the automorphism
group of the variety $X_{dI,t}$.
The coordinate functions of a Shioda map $\phi_{A,t}$ are products of the $y_i$ and hence if $\phi_A(y)=x$ then
$$
\phi_A(g(y))=(\zeta^{a'_1}x_1:\ldots:\zeta^{a'_n}x_n),\qquad
\mbox{where}\quad
\phi_A(y)=x=(x_1:\ldots:x_n).
$$
As $x_j=\prod_{k=1}^ny_k^{b_{jk}}$, the row vector $a'\in (\ZZ/d\ZZ)^n$ is obtained from the row vector $a\in (\ZZ/d\ZZ)^n$
as $a'=Ba$. Thus we get a homomorphism
$$
\Gamma_{d}\,\longrightarrow\,Aut(X_{A,t}),\qquad
g_a\longmapsto g_{Ba}.
$$
The kernel (image resp.) of this homomorphism will be denoted by
$\Gamma_A$ ($H_A$ resp.), so $H_A\cong\Gamma_d/\Gamma_A$.

In \cite{DGJ} some subgroups of $H_A$ are considered (called
discrete or scaling symmetries), but the next theorem shows that there is some advantage in considering the group $H_A$.

Two rational maps between algebraic varieties $f_i:X\rightarrow Y_i$, $i=1,2$
are said to be birationally equivalent if there is a Zariski open subset $U$ of $X$
and there are Zariski open subsets $U_i\subset Y_i$ with an isomorphism $\phi:U_1\rightarrow U_2$ such that $\phi\circ f_1=f_2$ on $U$.

\subsection{Theorem}\label{quotients}
Let $A$ be an $n\times n$ matrix with non-negative
integer entries such that the sum of the entries in any row and column
is equal to $n$ and such that $X_{A,t}$ is irreducible.
Then:

$\phi_{A,t}:X_{dI,t}\longrightarrow X_{A,t}$,
is birational to the quotient map
$X_{dI,t}\longrightarrow X_{dI,t}/\Gamma_A$,

$q_{A,t}:X_{A,t}\longrightarrow \overline{M}_{t}$,
is birational to the quotient map
$X_{A,t}\longrightarrow X_{A,t}/H_A$, and thus

$q_{A,t}\circ \phi_{A,t}:X_{dI,t}\longrightarrow\overline{M}_{t}$,
is birational to the quotient map
$X_{dI,t}\longrightarrow X_{dI,t}/\Gamma_d$.

\ts
In case $A=nI$, the Shioda map $\phi_{A,t}$
is given by $x_j=y_j^m$ and thus
$\Gamma_A=\mu_m^{n-1}$, the subgroup of $\Gamma_d$ generated
by $g_0,g_0g_1^n,\ldots,g_0g_1^n\cdots g_{n-2}^n$.
Similar to the argument in section \ref{mirrors}, the $\Gamma_A$-invariants (on each affine open subset ``$y_k=1$")
are generated by the $y_j^m$ and thus $\phi_{nI,t}$ is the quotient map.

The quotient group $H_A=\Gamma_d/\Gamma_A$ is now isomorphic to the group generated by $h_1=g_{(-m,m,0,\ldots)},\ldots, h_{n-2}=g_{(-m,0,\ldots,0,m,0)}$, these elements have order $d/m=n$.
Proceeding once more as in section \ref{mirrors} and comparing with \ref{gen mirror}, one finds that $q_{nI,t}$ is the quotient map and
$\overline{M}_t=X_{nI,t}/H_{nI}$.

Therefore the composition
$q_{nI,t}\circ \phi_{nI,t}$ is the quotient map
$X_{dI,t}\rightarrow \overline{M}_t$ and thus $\overline{M}_t=X_{dI,t}/\Gamma_d$.
Moreover, this composition
is given by
$$
z_0\,=\,x_1\cdots x_n\,=\,(y_1\cdots y_n)^m,\qquad
z_i\,=\,x_i^n\,=\,y_i^d.
$$

Now we verify that, for any $A$ as in the theorem, we have
$q_{A,t}\circ \phi_{A,t}=q_{nI,t}\circ \phi_{nI,t}$:
$$
z_0\,=\,x_1\cdots x_n\,=\,
\prod_{j=1}^n\left(\prod_{k=1}^ny_k^{b_{jk}}\right)\,=\,
\prod_{k=1}^ny_k^m,\qquad
z_i\,=\,\prod_{j=1}^5x_j^{a_{ij}}\,=\,
\prod_{k=1}^ny_k^{\sum_ja_{ij}b_{jk}}
\,=\,y_i^d.
$$

It remains to show that, for any $A$, the maps $\phi_{A,t}$ and $q_{A,t}$
are quotient maps. We do this by comparing degrees of maps (that is, the number of points in a general fiber). In case a map is a ``quotient by a finite group $G$ map", its degree is just the order $\sharp G$ of the group.

By definition of $\Gamma_A$ and the universal property of quotient varieties, there is a map
$X_{dI,t}/\Gamma_A\rightarrow X_{A,t}$ (basically: the coordinate functions
of the map $\phi_{A,t}$ are $\Gamma_A$-invariant and hence are functions on $X_{dI,t}/\Gamma_A$).
Thus $\mbox{deg}(\phi_{A,t})\geq \sharp \Gamma_A$.
The map $q_{A,t}:X_{A,t}\rightarrow \overline{M}_t$
factors over $X_{A,t}/H_A$ because the action of $H_A=\Gamma_d/\Gamma_A$ on $X_{A,t}$ is induced from the one of $\Gamma_d$ on $X_{A,t}$ and
$ \Gamma_d$ acts trivially on $\overline{M}_t=X_{dI,t}/\Gamma_d$.
Thus $\mbox{deg}(q_{A,t})\geq \sharp H_A$. Therefore
$$
\sharp\Gamma_d\,=\,\mbox{deg}(q_{nI,t}\circ \phi_{nI,t})\,=\,
\mbox{deg}(q_{A,t}\circ \phi_{A,t})\,=\,
\mbox{deg}(q_{A,t})\mbox{deg}(\phi_{A,t})\,\geq\,
(\sharp\Gamma_A)(\sharp H_A),
$$
and thus $\geq$ must be an equality, hence $\mbox{deg}(\phi_{A,t})= \sharp \Gamma_A$, $\mbox{deg}(q_{A,t})= \sharp H_A$ and the theorem follows.
\qed

\subsection{Differential forms}\label{diff forms}
The vector space of holomorphic $(n-1)$-forms on a smooth hypersurface $X=Z(F)$ of degree $d$ in $\PP^{n-1}$ has a basis
of the form
$$
\omega_{k,F}:=
\mbox{Res}_{X}\left(
y_1^{k_1}y_2^{k_2}\cdots y_n^{k_n}\frac{\sum_{i=1}^n(-1)^iy_i{\rm d}y_1\wedge\ldots\wedge
\widehat{{\rm d}y_i}\wedge\ldots\wedge{\rm d}y_n}{F}\right).
$$
with $k=(k_1,\ldots,k_n)$ and $k_i\in\ZZ_{\geq 0}$, $\sum k_i=d-n$.

Let $A$ be as Theorem \ref{quotients}.
As $X_{A,t}$ is CY, there is a unique holomorphic $(n-1)$-form, up to scalar multiple, we take
$$
\omega_{A,t}\,:=\,{\rm Res}_{X_{A,t}}(\tilde{\omega}_{A,t}),\qquad
\tilde{\omega}_{A,t}\,:=\,\omega_{0,F_{A,t}},
$$
so $k=0=(0,\ldots,0)$.
A somewhat tedious computation computes the pull-back of $\omega_{0,F_{A,t}}$ along the Shioda map $\phi_{A,t}:X_{dI,t}\rightarrow X_{A,t}$, the result implies that
$$
\phi_{A,t}^*\omega_{{A,t}}\,=\, c_A\omega_{(l,\ldots,l),F_{dI,t}},
\qquad l=m-1,\quad c_A=\det(B)/m.
$$
Note that this pull-back is essentially independent of $A$, in fact from
Theorem \ref{quotients} it follows that for any $A$ the pull-back must
be equal to the pull-back of the holomorphic $(n-1)$-form on $\overline{M}_t$
along $q_{A,t}\circ \phi_{A,t}=q_{nI,t}\circ \phi_{nI,t}$.

A more elegant way to obtain this differential form on $X_{dI,t}$ is to use the fact that it is, up to a scalar multiple, the only holomorphic
$(n-1)$-form on $X_{dI,t}$ which is invariant for the group $\Gamma_d$.
Using the description of the holomorphic $(n-1)$-forms given above, it is easy to see that indeed $\omega_{(l,\ldots,l),F_{dI,t}}$
spans the $\Gamma_d$-invariant forms.

\section{An example: the second family}

\subsection{The matrices $A$ and $B$}\label{AB}
Consider the Calabi-Yau variety defined by
$$
F_A :=
x_1^{n-1}x_2 + \ldots + x_i^{n-1}x_{i+1} + \ldots + x_n^{n-1}x_1,
$$
with the matrix $A$ given by:
$$
A = \begin{bmatrix}
	n-1 & 1 & 0 & \ldots & 0 \\
	0 & n-1 & 1 & \ldots & 0 \\
	0 & 0 & n-1 & \ldots & 0 \\
	\vdots & \vdots & \vdots & \ddots & \vdots \\
	1 & 0 & 0 & \ldots & n-1 \\
\end{bmatrix}.
$$
Note that sum of the entries in each row and each column is $n$.
Expanding by minors using the first column, there are two minors with non-zero determinant and these are upper diagonal, so one easily computes
$\det(A)=\left(n-1\right)^{n} - \left(-1\right)^{n}$.
The matrix $B$ is:
$$
B := \det\left(A\right) A^{-1} = \begin{bmatrix}
        q_1 & q_2 & q_3 & \ldots & q_n \\
	q_n & q_1 & q_2 & \ldots & q_{n-1} \\
	q_{n-1} & q_n & q_1 & \ldots & q_{n-2} \\
	\vdots & \vdots & \vdots & \ddots & \vdots \\
	q_2 & q_3 & q_4 & \ldots & q_1
\end{bmatrix},\qquad 
\mbox{with}\quad
q_i := \left(-1\right)^{i-1}\left(n-1\right)^{n-i}.
$$
Thus the Shioda map $\phi_{A,t}:X_{dI,t}\rightarrow X_{A,t}$
is given by:
$$
\phi_{A,t}\,:\,(y_1:\ldots:y_n)\,\longmapsto\,
(x_1:\ldots:x_n)\,=\,
(y_1^{q_1}y_2^{q_2}\cdots y_n^{q_n}:
 y_1^{q_n}y_2^{q_1}\cdots y_n^{q_{n-1}}:\ldots).
$$

As some of the exponents are negative, it is convenient to 
multiply all coordinate functions by a suitable (``minimal") monomial such that they become polynomials.
The smallest (negative) integer among the $q_i$ is $q_2=-(n-1)^{n-2}$,
thus we multiply by the monomial $(y_1y_2\cdots y_n)^{-q_2}$. As
$$
q_i-q_2=(-1)^{i-1}(n-1)^{n-i}+(n-1)^{n-2}\, \equiv \,
(-1)^{i-1+n-i}+(-1)^{n-2}\,
\equiv \,0 \, \,  \mbox{mod} \, \, n,
$$ 
this has the surprising side effect that the coordinate functions of
$\phi_{A,t}$ are polynomials in the $y_i^n$. Hence the Shioda map 
factors as $X_{dI,t}\rightarrow X_{mI, t}$, with $m=d/5$, given by
$y_i\mapsto y_i^n$ followed by a map $X_{mI, t}\rightarrow X_{A,t}$.
As the map $\phi_{5I,t}:X_{dI,t}\rightarrow X_{5I,t}$ obviously also
factors in this way, we get a diagram:
$$
\begin{array}{rlcrl}
&&X_{dI,t}&&\\
&\!\swarrow &  \!\downarrow &  \searrow\! & \\
\phantom{|} X_{5I,t} \phantom{|} & \leftarrow & \phantom{|} X_{mI, t}\phantom{|} & \rightarrow & \phantom{|} X_{A,t} \phantom{|}, \\
\end{array} 
$$
where the vertical map is given by $u_k = y_k^n$ and $u_k$ are coordinates on $X_{mI,t}$,
note that the relation $F_{dI,t}(y_1,\ldots,y_5)=0$ implies
$F_{mI,t}(u_1,\ldots,u_5)=0$.
Finally, notice that an argument, 
similar to the one which showed that 
$X_{dI,t} \rightarrow X_{A,t}$ is a quotient map, 
shows that all maps in the diagram above are quotient maps. 
Denote by $\mu_A$ the group corresponding to the vertical map and by $\Gamma_A'$ the group corresponding to the lower right map.

For example, in case $n=5$, we have $d=1025$, $q_2=-64$, and, with $u_i=y_i^5$,
{\renewcommand{\arraystretch}{1.4}
$$
\begin{array}{rcl}
\phi_{A,t}\,:\,(y_1:\ldots:y_n)\,\longmapsto\,
(x_1:\ldots:x_n)&=&
(y_1^{256}y_2^{-64}y_3^{16}y_4^{-4}y_5:
 y_1^{}y_2^{256}y_3^{-64}y_4^{16}y_5^{-64}:\ldots)\\
&=&(y_1^{320}y_2^{0}y_3^{80}y_4^{60}y_5^{65}:
 y_1^{65}y_2^{320}y_3^{0}y_4^{80} y_5^{60}:\ldots)\\
&=&
(u_1^{64}u_3^{16}u_4^{12}u_5^{13}:
 u_1^{13}u_2^{64}u_4^{16}u_5^{12}:\ldots).
\end{array}
$$
}

\subsection{The groups $\Gamma_A$ and $H_A$ in case $n=5$}
The group $H_A$ is the image of $\Gamma_d$ under the homomorphism $g_a\mapsto g_{Ba}$. Denoting by $\hat{g}_i$ the image of the generator
$g_i$ of $\Gamma_d$ given in section \ref{autos}, one easily verifies that (note that the entries are taken mod $d=1025$):
$$
\hat{g}_0\,=\,g_{(1280,5,-20,80,-320)}\,=\,g_{(255,5,1005,80,705)}.
$$
As $g_0^{41}=g_{(205,205,205,205,205)}$,
the automorphism $\hat{g}_0$ has order $41$, and one also finds:
$$
\hat{g}_1\,=\,\hat{g}_0^{51}\;(=\hat{g}_0^{10}),\qquad
\hat{g}_2\,=\,\hat{g}_0^{-13},\qquad
\hat{g}_3\,=\,\hat{g}_0^{3}.
$$
As $\Gamma_d\cong(\ZZ/205)\times(\ZZ/1025)^3$ we get
$$
H_A\,=\,\langle \hat{g}_0\rangle\,\cong\, \ZZ/41\ZZ,\qquad \Gamma_A\,\cong\,(\ZZ/5\ZZ)\times(\ZZ/1025\ZZ)^3.
$$
Furthermore, it is easy to check that, with $\mu_A$ and $\Gamma'_A$
as in section \ref{AB},
$$
\mu_A \,\cong\, (\ZZ/5\ZZ)^3, \qquad 
\Gamma_A' \,\cong\, \ZZ/5\ZZ \times (\ZZ/205\ZZ)^3.
$$

In \cite{DGJ}, table 1, an automorphism of order $41$ of $X_{A,t}$ is given.
In our notation, based on $d$-th roots of unity with
$d=1025=25\cdot 41$, it is $g_b$ with $b:=25\cdot(1,37,16,18,10)$.
One verifies that  $b\equiv 185\cdot(255,5,1005,80,705)$ mod $1025$. Thus
$\hat{g}_0^{185}$ is the automorphism of $X_{A,t}$ considered in \cite{DGJ}. As a consequence of Theorem \ref{quotients} we obtain
that $X_{A,t}/H_A$, with $H_A\cong\ZZ/41\ZZ$, is birationally isomorphic to $\overline{M}_t$.
This was already suggested in \cite{DGJ} (1.8).

\subsection{The Shioda correspondence}\label{correspondence}
In the papers \cite{GPR} and \cite{DGJ}, the Dwork pencil $X_{5I,t}$
and the pencil $X_{A,t}$ are related by fractional changes of variables. From a geometrical point of view, this is just the correspondence between these pencils defined by the pencil $X_{dI,t}$. Indeed, using coordinates $x'_j$ for the Dwork pencil, the map
$\phi_{5I,t}:X_{dI,t}\rightarrow X_{5I,t}$ is given by
$x'_j=y_j^{205}$, so $y_j=(x'_j)^{1/205}$. The map
$\phi_{A,t}:X_{dI,t}\rightarrow X_{A,t}$ is given by
$x_j=\prod y_k^{b_{jk}}=\prod (x'_k)^{b_{jk}/205}$.
Thus if we redefine (as in \cite{DGJ}) $y_j:=x_j$ and $x_j:=x_j'$
and use that $b_{1k}=q_k$ then
(with a shift in indices in our notation w.r.t.\ the one in \cite{DGJ}, (1.7)) we find the ``change of variables" from \cite{DGJ}, (1.7):
$$
y_1=x_1^{\frac{256}{205}}x_2^{\frac{-64}{205}}x_3^{\frac{16}{205}}
x_4^\frac{-4}{205}x_5^\frac{1}{205},\qquad\mbox{etc.}
$$

Thus the Shioda correspondence below and the induced map on cohomology $\psi$ (note that since $\phi_{A,t}$ is in general only a rational map one has to do some blowing up on $X_{dI,t}$),
{\renewcommand{\arraystretch}{1.3}
$$
\begin{array}{rlcrl}
&&X_{dI,t}&&\\
\phantom{|}^{\phi_{5I,t}}\phantom{|}&\!\swarrow&&\searrow\!&
\phantom{|}^{\phi_{A,t}}\phantom{|}\\
X_t\;&&&&\;X_{A,t},\\
\end{array}\qquad
\psi\,:=\,(\phi_{A,t})_*(\phi_{5I,t})^*\,:\,
H^3(X_{5I,t},\QQ)\,\longrightarrow\, H^3(X_{A,t},\QQ)
$$
}

\noindent
seem to play an important role in Mirror symmetry.
The linear map $\psi$
is a map of rational Hodge structures and it
induces an isomorphism on $H^{3,0}$ (cf.\ \ref{diff forms}).
Hence it gives another explanation (besides the one offered in \ref{PF})
for the equality of the Picard Fuchs equations for the holomorphic
three forms in both pencils.

\end{document}